\documentclass{amsart}

\usepackage{amsfonts, amsmath, amssymb}

\newtheorem{prop}{Proposition}[section]
\newtheorem{theo}[prop]{Theorem}
\newtheorem{coro}[prop]{Corollary}
\newtheorem{lem}[prop]{Lemma}

\theoremstyle{definition}
\newtheorem{dfn}[prop]{Definition}
\newtheorem{exam}[prop]{Example}

\theoremstyle{remark}
\newtheorem{rem}[prop]{Remark}

\def\<{\langle}
\def\>{\rangle}

\def\l{{\it l}}
\def\s{{\sigma}}
\def\Si{{\Sigma}}
\def\G{{\Gamma}}
\def\t{{\tau}}

\def\O{{\mathcal O}}
\def\P{{\mathbf P}}
\def\R{{\mathbf R}}

\def\C{{\mathbf C}}
\def\Z{{\mathbf Z}}
\def\Q{{\mathbf Q}}
\def\L{{\mathcal L}}
\def\K{{\mathcal K}}
\def\F{{\mathcal F}}
\def\E{{\mathcal E}}
\def\logD{{\rm log}\, D}

\begin{document}

\title{The Bott Formula for Toric Varieties}

\author{Evgeny~N.~Materov}

\date{This edition: \today}

\address{Krasnoyarsk State Technical University, 
Sverdlova 5, Zheleznogorsk 662990, Russia}

\email{materov@atomlink.ru}

\curraddr{Mathematisches Institut, Universit\"at T\"ubingen,
Auf der Morgenstelle 10, T\"ubingen 72076, Germany}

\email{evgeny.materov@uni-tuebingen.de}

\thanks{2000 {\em Mathematics Classification.} Primary 14M25, 32L10, 58A10; 
Secondary 52B20, 52B11.
\newline
Supported by RFFI, grant 00-15-96140 and by DFG, 
Forschungsschwerpunkt "Globale Methoden der komplexen Geometrie.}





\begin{abstract}
The purpose of this paper is to give an explicit formula which allows one to
compute the dimension of the cohomology groups of the sheaf
$\Omega_{\P}^p(D)=\Omega_{\P}^p\otimes {\O_\P}(D)$
of $p$-th differential forms of Zariski twisted by an ample invertible sheaf on
a complete simplicial toric variety. The formula involves some combinatorial
sums of integer points over all faces of the support polytope for
${\O_\P}(D)$.
Also, we introduce a new combinatorial object, the so-called
$p$-th Hilbert-Ehrhart polynomial, which generalizes the usual notion and
behaves similar. Namely, there exists a generalization of the reciprocity law
for a usual Hilbert-Ehrhart polynomial.
Some applications of the Bott formula are discussed.
\end{abstract}


\maketitle


\section{Introduction}

One of the most important invariants in algebraic geometry deals with
the computation of cohomologies of sheaves.
The crucial role in this computation is played by some vanishing criteria for
higher-dimensional cohomologies of certain sheaves.
This vanishing then enables one to reduce the calculation of the whole
cohomology group $H^q(\P,\F)$ to the vector space of global sections
$\G(\P,\F)=H^0(\P,\F)$. For example, the Bott vanishing theorem \cite{Bott}
and the Serre duality theorem help to compute the dimension of
$H^q(\P^n,\Omega_{\P^n}^p(k))$ on a complex projective space
$\P^n=\P^n(\C)$, where
$\Omega_{\P^n}^p(k):=\Omega_{\P^n}^p\otimes_{\O_{\P^n}} \O_{\P^n}(k)$.

For simplicity, we denote the dimension of the $q$-th cohomology group of the
coherent sheaf $\F$ on a complete projective variety $\P$ by
$h^q(\P,\F) := \dim_{\C} H^q(\P,\F)$.
\begin{theo}[Bott formula for $\P^n$ {\rm \cite{Bott, Okonek}}]
Let $\P^n$ be a complex projective space.
\begin{itemize}
 \item[{\rm 1)}] If $k = 0$, then
   \begin{eqnarray*}
      h^q(\P^n,\Omega_{\P^n}^p)=
      \left\{
      \begin{array}{ll}
       \displaystyle 
       1,& p=q,\\
       \\
       0,&otherwise.
      \end{array}
      \right.
     \end{eqnarray*}
 \item[{\rm 2)}] If $k \ne 0$, then
    \begin{eqnarray*}
      h^q(\P^n,\Omega_{\P^n}^p(k))=
      \left\{
      \begin{array}{ll}
       \displaystyle
       \binom{n+k-p}{k}\binom{k-1}{p},&q=0, \:k>p,\\
       \\
       \displaystyle
       \binom{-k+p}{-k}\binom{-k-1}{n-p},&q=n,\: k<p-n,\\
       \\
       0,&otherwise.
      \end{array}
      \right.
    \end{eqnarray*}
\end{itemize}
\label{Original.Bott}
\end{theo}
This theorem was also proved by I-C.~Huang \cite{Huang}, who constructed an 
explicit basis for the cohomology groups. 
Theorem 2.3.2 in \cite{Dolgachev} gives a 
similar result for weighted projective spaces.

The first aim of this paper is to give a generalization of the Bott formula
for a complete projective toric variety.
Let $D$ be a Cartier divisor on a toric variety $\P$, and 
$\Omega_\P^p(D)$ be the sheaf of $p$-differential
forms of Zariski on $\P$.
Since we have an analogy of the Bott vanishing
theorem, it is possible to compute the dimension of
$H^q(\P,\Omega_\P^p(D))$ when $D$ is ample.
The answer will be given in terms of combinatorial sums over all faces
of the {\it support polytope} $\Delta$ for $\O_{\P}(D)$, playing the 
role of the number $k$.
The case when $\O_{\P}(D)\simeq\O_{\P}$ corresponding to $k=0$ was
extended to a toric variety in \cite{Danilov-Khovanskii} and
\cite{Oda.book}.
Also, note that for $p=0$, we get the well-known result on the cohomologies
$H^q(\P,\O_{\P}(D))$ of the ample invertible sheaf $\O_{\P}(D)$, 
which can be found in any introduction to the theory of toric varieties:
\[
  \dim_\C H^q(\P,\O_{\P}(D))=
  \left\{
  \begin{array}{ll}
    \# (\Delta\cap M),&q=0,\\
    \\
    0                ,&q>0,
  \end{array}
  \right.
\]
where $\# (\Delta\cap M)$ is the number of integer points in the polytope 
$\Delta$ in a lattice $M$.

Let $N$ be a free $\Z$-module of rank $n$ and $M := {\rm Hom}_\Z(N,\Z)$ its
dual. Denote by $\Si$ a complete {\it fan} of convex polyhedral cones in
$N_\R := N\otimes \R$. We associate a so-called {\it toric variety}
$\P = \P(\Si)$ with $\Si$, i.e. a normal complex variety containing a torus
$T_N := {\rm Hom}_\Z(M,\C^*)$ as a dense open set with an algebraic action
of $T_N$ on $\P$. For precise definitions and basic references about toric
varieties see \cite{Danilov}, \cite{Oda.book} and \cite{Fulton.toric}.
In section \ref{section.cohomologies} we will give a brief
review of another, but equivalent, approach to toric varieties which
has been introduced and studied in \cite{Audin}, \cite{Cox} and
\cite{Batyrev-Cox}. Throughout this paper we assume that $\P$ is complete
and simplicial. Here is the precise definition of the sheaf of $p$-differential
forms of Zariski, which is the main object of our study.
\begin{dfn}
Denote {\it the sheaf of $p$-differential forms of Zariski
on $\P$} by $\Omega_{\P}^p$.
These forms are defined by $\Omega_{\P}^p:=i_{*}\Omega_U^p$, where
$i:\,U\hookrightarrow \P$ is the inclusion of the nonsingular locus 
$U$ of $\P$.
Also, denote for any invertible sheaf $\O_{\P}(D)$ on $\P$
\[
  \Omega_{\P}^p(D):=\Omega_{\P}^p\otimes \O_{\P}(D).
\]
\end{dfn}

The paper is organized as follows:

{\it In sections 2 and 3}, we prove a generalization of the original 
Bott formula stated above.
The generalization will be given in two forms:
in theorem \ref{Bott.formula.I} and in theorem
\ref{Bott.formula.II}. The first version of the Bott formula is proved
using the Ishida-Oda complex, while the second version is  
proved using the technique of Danilov and Khovanski\^i.
Both versions are formulated in terms of sums over all faces $\Gamma$ of the 
support polytope $\Delta$.

{\it In section 4}, we compare our results with the original Bott formula.

{\it In section 5}, we introduce the notion of the $p$-th Hilbert-Ehrhart 
polynomial $L_p(k)$ coinciding with the usual Ehrhart polynomial $L(k)$
when $p=0$. Our idea is to compare the Ehrhart polynomial with a 
certain Hilbert polynomial corresponding to the sheaf $\Omega_{\P}^p(D)$.
Then, the Serre duality provides a generalization of the well-known 
{\it reciprocity law}.

{\it In section 6}, we compare the two versions of the Bott formula for 
toric varieties and obtain non-trivial identities between integer points
in simple integer polytopes. We give a simple combinatorial proof
of the identities and the generalized reciprocity law.

{\it In section 7}, we apply the Bott formula to the computation of
cohomologies of quasi-smooth hypersurfaces on a complete simplicial 
toric variety.

{\it In section 8}, we compute the dimension of the space of  
global sections of weight components of the sheaf 
$\Omega_{\P}^p({\rm log}(-K))\otimes \L$.

{\it In section 9}, we obtain a result similar to the Bott formula 
on the projective bundle $\P(\E)$ associated with the 
sheaf $\E = \L_0\oplus \ldots \oplus \L_s$, where $\L_i$ is an ample
invertible sheaf on $\P$. Namely, we compute the dimension of the 
global sections of the sheaf $\Omega_{Y/\P}^p(k)$ of relative differential 
forms on $\P(\E)$.

{\it Acknowledgments.} The author would like to thank Prof. A.K.~Tsikh for his
wise advices and very useful discussions during the preparation of this paper.
I am very grateful to Prof. A.P.~Yuzhakov, who introduced me to the theory of
combinatorial sums. I appreciate the hospitality of the Mathematisches 
Institut of Universit\"at T\"ubingen, where this paper was finished.

\section{The generalized Bott formula (I)}


Let $\P$ be a complete simplicial $n$-dimensional projective toric variety.
Fix an invertible sheaf $\L=\O_{\P}(D)$ corresponding to an ample Cartier 
divisor $D$ on $\P$.

It is known from the general theory of toric varieties
(see, e. g., \cite{Oda.book}) that the
cohomology group $H^0(\P,\L)$ admits an eigenspace decomposition
coming from the torus action
\[
  H^0 (\P,\L) = \bigoplus_{m\in M} H^0 (\P,\L)_m
\]
and $H^q(\P,\L)$ vanishes for $q>0$.
\begin{dfn} Denote the set of all $k$-dimensional cones of $\Si$ by
$\Si(k)$ and by $|\Si(k)|$ the number of $k$-dimensional cones 
of $\Si$, i.e. the cardinality of $\Si(k)$.
\end{dfn}
\begin{dfn}[\cite{Torus.Embeddings}]
The convex hull $\Delta = \Delta(\L)$ of all lattice points $m\in M$ 
in $M_\R:=M\otimes \R$ for which $H^0 (X,\L)_m \ne 0$ is called the
{\it support polytope for $\L$}.
Moreover, $\Delta$ is equal to the intersection
\[
  \Delta = \bigcap_{\t\in\Si(n)} (m_\t(\L)+\t^\vee),
\]
where $m_\t(\L)$ is the unique lattice point of $M$ such that
\[
  \{ m\in M:\, H^0(U_\t,\L)_m \ne 0 \} =
     m_\t(\L) + \t^\vee \cap M
\]
for any affine chart $U_\t={\rm Spec}\C[M\cap \t^\vee]$
corresponding to the cone $\t$.
\end{dfn}
\begin{rem} $\Delta$ is a {\it simple} convex polytope of 
dimension $n$, i.e. only $n$ edges of $\Delta$ meet at each vertex, 
if and only if $\P$ is a simplicial toric variety.
\end{rem}
\begin{dfn}
Each face $\Delta_\s$ of $\Delta$ corresponding to the cone 
$\s\in\Si$ is defined as
\[
  \Delta_\s:=
  \bigcap_{\begin{subarray}{l}
             \t\in\Si(n) \\
             \t\succ\s
           \end{subarray}}
  (m_\t(\L) + \t^\vee \cap \s^\perp),
\]
where $\t\succ\s$ means that $\s$ is a face of $\t$.
\end{dfn}
\begin{dfn} Denote by $\l (\Delta)$ the number of integer points in the
polytope $\Delta$, and let $\l (\Delta_\s)$ be the number of integer points
contained in the face $\Delta_\s$ of $\Delta$.
\label{Int.Points}
\end{dfn}
\begin{dfn}
We identify the closure $V(\s)$ of any torus-invariant orbit 
$orb(\s)$ in $\P$ corresponding to $\s\in\Si$ with
a toric variety with respect to a fan glued from the cones
$(\t + (-\s))/\R\s$ in $N_\R/\R\t$, $\t\in {\rm Star}_\s(\Si)$,
where $\R\s:=\s+(-\s)$ is the smallest $\R$-subspace containing $\s$ of
$N_\R$, while ${\rm Star}_\s(\Si) := \{\t\in\Si:\,\t\succ\s\}$ 
(see \cite{Oda.book}).
\end{dfn}
\begin{prop}[\cite{Batyrev-Cox}, Proposition 4.10]
If $\L$ is an ample invertible sheaf on $\P$, then one has a one-to-one
correspondence between $(n-k)$-dimensional faces $\Delta_\s$ of the 
polytope $\Delta(\L)$ and $k$-dimensional cones $\s\in\Si$ reversing the 
face-relation. Moreover, $\Delta_\s$ is the support polytope for the 
sheaf $\O_{V(\s)}\otimes \L$.
\label{support.polytope}
\end{prop}
\begin{dfn}
We denote by $\F_k(\Delta)$ the set of all $k$-dimensional faces of $\Delta$
and by $f_k = |\F_k(\Delta)|$ the number of $k$-dimensional faces of $\Delta$.
\end{dfn}
\begin{rem} 
It follows from the proposition above that 
$\F_j(\Delta) \simeq \Si(n - j)$ and $f_j = |\Si(n - j)|$ respectively. 
\end{rem}
\begin{dfn} 
Define the torus-invariant effective divisor $D(\s)$ on each $V(\s)$
by
\[
  D(\s):=
  V(\s)\backslash orb(\s).
\]
\end{dfn}

Now recall the construction of the complex which has been introduced and 
studied by M.-N.~Ishida and T.~Oda in \cite{Ishida}, \cite{Oda.book} and
\cite{Oda.de.Rham}, known as {\it Ishida's $p$-th complex of
$\O_{\P}$-modules}. Ishida's complex plays an important role in the
proof of the first version of the generalized Bott formula.
\begin{prop}[\cite{Oda.book}, Corollary 3.2)]
Let $\Omega_{V(\s)}^p(\logD(\s))$ be a sheaf of meromorphic $p$-forms with 
logarithmic poles along $D(\s)$.
There exists an isomorphism of $\O_{V(\s)}$-modules
\[
  \Omega_{V(\s)}^p(\logD(\s))=
  \O_{V(\s)}\otimes_\Z\bigwedge\nolimits^p
  (M\cap \s^\perp) \qquad\mbox{for}\quad 0\le p\le \dim V(\s).
\]
\end{prop}
Now, for each pair of integers $p$ and $q$ we let
\[
  \K_{\P}^{p,\,q}:=
  \bigoplus_{\s\in\Si(q)}
  \Omega_{V(\s)}^{p-q}(\logD(\s))=
  \bigoplus_{\s\in\Si(q)}
  \O_{V(\s)}\otimes_\Z\bigwedge\nolimits^{p-q}
  (M\cap \s^\perp) \qquad\mbox{if}\quad 0\le q\le p,
\]
and $\K_{\P}^{p,\,q} = 0$ otherwise.
\begin{dfn} 
The coboundary map 
$
\delta:
\K_{\P}^{p,\,q}
\rightarrow
\K_{\P}^{p,\,q+1}
$
is defined to be the direct sum $\delta=\oplus_{\s,\t} R_{\s/\t}$
of the maps
\[
  R_{\s/\t}:\,
  \Omega_{V(\s)}^{p-q}(\logD(\s))
  \rightarrow
  \Omega_{V(\t)}^{p-q-1}(\logD(\t)),
\]
i.e.
\[
  R_{\s/\t}:\,
  \O_{V(\s)}\otimes_\Z\bigwedge\nolimits^{p-q}
  (M\cap \s^\perp)
  \rightarrow
  \O_{V(\t)}\otimes_\Z\bigwedge\nolimits^{p-q-1}
  (M\cap \t^\perp),
\]
where $R_{\s/\t}$ is the tensor product of the restriction map
$\O_{V(\s)}\rightarrow \O_{V(\t)}$ with the interior product $\delta_{\s/\t}$
whenever $\s$ is a face of $\t$, and zero otherwise. The definition of
$\delta_{\s/\t}$ is the following. The homomorphism $\delta_{\s/\t}$ is defined
to be zero when $\s$ is not a face of $\t$. On the other hand, if $\s$ is
a face of $\t$, then we can determine a primitive element $n\in N$
uniquely modulo $N\cap \R\s$ so that $\t+(-\s)=\R_{\ge 0}n+\R\s$.
Moreover, $M\cap \t^\perp$ is a $\Z$-submodule of corank one in
$M\cap \s^\perp$. Each element of
$\bigwedge\nolimits^{p-q}(M\cap \s^\perp)$ can be written as a finite linear
combination of
\[
  m_1\wedge m_2\wedge \cdots\wedge m_{p-q},
  \quad m_1\in M\cap \s^\perp,
  \quad m_2,\ldots, m_{p-q}\in M\cap \t^\perp.
\]
Now define
\[
  \delta_{\s/\t}:
  \bigwedge\nolimits^{p-q}   (M\cap \s^\perp)\rightarrow
  \bigwedge\nolimits^{p-q-1} (M\cap \t^\perp)
\]
by
\[
  \delta_{\s/\t}(m_1\wedge m_2\wedge \cdots\wedge m_{p-q}):=
  \<m_1,n\>m_2\wedge \cdots\wedge m_{p-q}.
\]
Actually, $R_{\s/\t}$ is the {\it Poincar\'e residue map} for the component
$V(\t)$ of the divisor $D(\s)$ on $V(\s)$.
\end{dfn}
\begin{theo}[\cite{Oda.book}, Theorem 3.6]
If $\P$ is a simplicial toric variety, then for each $0\le p\le n$ 
there exists an exact sequence
  \begin{equation}
    0\rightarrow \Omega_{\P}^p
     \rightarrow
     \K_{\P}^{p,\,0}
     \stackrel{\delta}{\rightarrow}
     \K_{\P}^{p,\,1}
     \stackrel{\delta}{\rightarrow}
     \ldots
     \stackrel{\delta}{\rightarrow}
     \K_{\P}^{p,\,p}
     \rightarrow 0
    \label{Ishida's.complex}
   \end{equation}
of $\O_{\P}$-modules on $\P$.
\end{theo}
\begin{theo}[The generalized Bott formula (I)] 
Let $\P$ be a complete simplicial toric variety, 
$\O_{\P}(D)$ be an invertible sheaf on $\P$ corresponding to an ample
divisor $D$, and $\Delta$ be the support polytope for $\O_{\P}(D)$. Then
\begin{itemize}
\item[{\rm 1)}] 
  \[
    h^q(\P,\Omega_{\P}^p) =
    \left\{
    \begin{array}{ll}
       \displaystyle
       \sum_{j=0}^p(-1)^{p-j}\binom{n-j}{p-j} f_{n - j},  &q=p,\\
       \\
       0,                                                 &q\ne p.
    \end{array}
    \right.
\]
\item[{\rm 2)}] 
  \[
    h^q(\P,\Omega_{\P}^p(D)) =
    \left\{
    \begin{array}{ll}
       \displaystyle
       \sum_{j=0}^p(-1)^j\binom{n-j}{p-j} 
       \sum_{F\in\F_{n - j}(\Delta)}\l(F),&q=0,\\
       \\
       0,                                 &q>0.
    \end{array}
    \right.
  \]
\end{itemize}
\label{Bott.formula.I}
\end{theo}
\begin{proof} The assertion $1)$ is proved in the Theorem 3.11 
of \cite{Oda.book} in the form 
  \[
    h^p(\P,\Omega_{\P}^p) =
    \displaystyle
    \sum_{j=0}^p(-1)^{p-j}\binom{n-j}{p-j} |\Si(j)|,
\]
and $h^q(\P,\Omega_{\P}^p) = 0$ for $p\ne q$.

Let us prove $2)$. We write $\L$ for $\O_{\P}(D)$,
and $\Omega_{\P}^p\otimes\L$ for $\Omega_{\P}^p(D)$ for convenience.
By Bott vanishing theorem (see \cite{Danilov, Oda.book})
$h^q(\P,\Omega_{\P}^p\otimes\L) = 0$ for $q > 0$. The sequence 
(\ref{Ishida's.complex}) remains exact after shifting by $\L$
\begin{equation}
\label{twisted.complex}
  0\rightarrow \Omega_{\P}^p\otimes \L
   \rightarrow
   \K_{\P}^{p,\,0}\otimes \L
   \rightarrow
   \K_{\P}^{p,\,1}\otimes \L
   \rightarrow
   \ldots
   \rightarrow
   \K_{\P}^{p,\,p}\otimes \L
   \rightarrow 0.
\end{equation}
Note, that $M\cap \s^\perp$ is a free $\Z$-module of rank $n-j$ and
$\bigwedge\nolimits^{p-j}(M\cap \s^\perp)$ is a free $\Z$-module of
rank $\displaystyle \binom{n-j}{p-j}$ when $\s\in\Si(j)$.
This observation shows that
\[
  \K_{\P}^{p,\,j}\otimes \L
  \simeq
  \O_{V(\s)}       \otimes_\Z\bigwedge\nolimits^{p-j}
  (M\cap \s^\perp)\otimes \L
  \simeq
  \displaystyle
  \O_{V(\s)}\otimes \L^{\oplus\binom{n-j}{p-j}}.
\] 
Assume for the moment that $h^q(\P,\K_{\P}^{p,\,j}\otimes \L) = 0$ for 
$q > 0$. Then (\ref{twisted.complex}) gives the exact sequence of 
global sections
\begin{eqnarray*}
  0\rightarrow H^0(\P,\Omega_{\P}^p\otimes \L)
   \rightarrow H^0(\P,\K_{\P}^{p,\,0}\otimes \L)
   \rightarrow H^0(\P,\K_{\P}^{p,\,1}\otimes \L)
   \rightarrow
   \ldots\\
   \ldots
   \rightarrow H^0(\P,\K_{\P}^{p,\,p}\otimes \L)
   \rightarrow 0.
\end{eqnarray*}
Since $h^0(\P,\O_{V(\s)}\otimes \L)=l(\Delta_{\s})$ 
(see Proposition \ref{support.polytope}), we have
 \begin{eqnarray*}
   h^0(\P,\Omega_{\P}^p\otimes \L) 
   &=&
   \sum_{j=0}^p (-1)^j h^0(\P,\K_{\P}^{p,\,j}\otimes \L)\\
   &=&
   \sum_{j=0}^p (-1)^j \binom{n-j}{p-j}
   \sum_{\s\in\Si(j)}
   h^0(\P,\O_{V(\s)}\otimes \L)\\
   &=&
   \sum_{j=0}^p(-1)^j\binom{n-j}{p-j} 
   \sum_{\s\in\Si(j)}\l(\Delta_\s)\\
   &=&
   \sum_{j=0}^p(-1)^j\binom{n-j}{p-j} 
   \sum_{F\in\F_{n - j}(\Delta)}\l(F).
 \end{eqnarray*}

It remains to show that $h^q(\P,\K_{\P}^{p,\,j}\otimes \L) = 0$ for 
$q > 0$. Moreover, it is sufficient to prove the vanishing of 
$h^q(\P,\O_{V(\s)}\otimes \L) = h^q(\P,\L|_{V(\s)})$. 
As observed in \cite{Danilov}, Lemma 6.8.1, the restriction homomorphism
\[
 \Gamma(\P,\L)\rightarrow \Gamma(\P,\L|_{V(\s)})
\]
is surjective. Hence, the sheaf $\L|_{V(\s)}$ is generated by its 
global sections since $\L$ is ample, and all cohomologies 
$h^q(\P,\L|_{V(\s)})$ for $q > 0$ vanish by the Bott vanishing theorem. 
\end{proof}
\begin{coro} One has the following formulas
\[
  \begin{array}{lcl}
  &&\displaystyle
  \sum_{p=0}^n h^p(\P,\Omega_{\P}^p)\, y^p =
  \displaystyle
  \sum_{j=0}^n f_j (1 - y)^j y^{n - j}=
  \displaystyle
  \sum_{j=0}^n f_j (y - 1)^j,\\
  &&\displaystyle
  \sum_{p=0}^n h^0(\P,\Omega_{\P}^p(D))\, y^p =
  \displaystyle
  \sum_{j=0}^n \Delta^{(j)}
       (- y)^j (y + 1)^{n - j},
  \quad \mbox{\it if}\; D\; \mbox{\it is ample},
\end{array}
\]
where $\displaystyle \Delta^{(j)} := \sum_{F\in\F_{n - j}(\Delta)}\l(F)$.
\label{Bott.Corollary.I}
\end{coro}
Notice that the first line of the equalities is contained in  
Theorem~3.11 of \cite{Oda.book}, where the second equality follows from the
Serre duality $h^p(\P,\Omega_{\P}^p)=h^{n-p}(\P,\Omega_{\P}^{n-p})$.

\section{The generalized Bott formula (II)}

Another approach gives the second version of the generalization of the
Bott formula. Here we present two independent proves of this formula.

Suppose that $\P$ is a complete simplicial $n$-dimensional toric variety 
as above.
Assume that $\O_{\P} (D)$ is an ample invertible sheaf on $\P$ determining
the convex polytope $\Delta$.
Because of the one-to-one correspondence between the faces 
$\Gamma = \Delta_\s$ of the polytope $\Delta$ and cones $\s$ of the fan $\Si$,
we can associate a "small" toric variety, i.e. the closure of torus-invariant 
orbit of $\P$, with any face $\Gamma\subseteq\Delta$.
\begin{dfn} 
The toric $k$-dimensional subvariety $\P_{\Gamma}$ of 
$\P$ corresponding to the $k$-dimen\-sional face $\Gamma\subseteq\Delta$
is the closure of a $n - k$-dimensional orbit $orb(\Gamma)$ in $\P$.
For $\Gamma = \Delta$ we let $\P_{\Delta} = \P$.
Define the effective divisor $D_{\Gamma}$ on each $\P_{\Gamma}$ by
\[
  D_{\G} := \P_{\Gamma}\backslash orb(\G).
\]
\end{dfn}
\begin{dfn} 
\cite{Danilov-Khovanskii, Danilov2}. Denote by
$\displaystyle\Omega_{(\P_{\Gamma},\,D_{\Gamma})}^p$ the sheaf of regular
differential $p$-forms on the toric variety $\P_{\Gamma}$ with zeros along
$D_{\Gamma}$, arising from the exact sequence
\begin{equation}
  0\rightarrow
   \Omega_{(\P_{\Gamma},\,D_{\Gamma})}^p
   \rightarrow
   \Omega_{\P_{\Gamma}}^p
   \stackrel{R_\Gamma}{\rightarrow}
   \bigoplus_{\Theta}\Omega_{\P_{\,\Theta}}^p
   \rightarrow 0,
\label{restr.hom}
\end{equation}
where $\Theta$ runs over all facets (faces of codimension one) of 
$\Gamma$ and $R_{\Gamma}$ is the restriction homomorphism.
\end{dfn}
\begin{dfn} 
The {\it Euler-Poincar\'e characteristic}
of a coherent $\O_{\P}$-module $\F$ on a $n$-dimensional complete
variety $\P$ over $\C$ is defined to be
\[
  \chi(\P,\F) :=
  \sum_{k = 0}^n (-1)^k \dim_{\C} H^k (\P,\F)=
  \sum_{k = 0}^n (-1)^k h^k (\P,\F).
\]
\end{dfn}
\begin{dfn}
Denote by $\l^{*}(\Delta)$ the number of integer points in
the relative interior of the polytope $\Delta$.
Also, let $\l^{*}(\Gamma)$ be the number of integer points
contained in the relative interior of a face $\Gamma$ of $\Delta$.
\end{dfn}

Notice the following key property of the sheaf
$\Omega_{(\P_{\Gamma},\,D_{\Gamma})}^p$.
\begin{prop} Let $\P$ be a complete simplicial toric variety and 
$\P_{\Gamma}$ be a closed subset of $\P$ associated with the face 
$\Gamma$ of the support polytope $\Delta$ corresponding to an 
invertible sheaf $\O_{\P}(D)$ on $\P$.
Then
\begin{eqnarray*}
  \chi(\P_{\Gamma},\,\Omega_{(\P_{\Gamma},\,D_{\Gamma})}^p
  \otimes \O_{\P}(D)) =
  \left\{
  \begin{array}{ll}
     \displaystyle
     (-1)^{\dim\Gamma}\binom{\dim\Gamma}{p}, &\O_{\P}(D)\simeq \O_{\P},\\
     \\
     \displaystyle
     \binom{\dim\Gamma}{p} \l^{*}(\G),       &\O_{\P}(D)\not\simeq \O_{\P}, 
                                             \; D \; \mbox{\it is ample}.\\
  \end{array}
  \right.
\end{eqnarray*}
\label{Remark.in.DK}
\end{prop}
\begin{proof} The statement follows immediately from the Proposition
2.10 of \cite{Danilov-Khovanskii} by taking the restriction of the sheaf 
$\L$ to $\P_{\Gamma}$.
\end{proof}

Here is the main result of this section.
\begin{theo}[The generalized Bott formula (II)] 
Let $\P$ be a complete simplicial toric variety, 
$\O_{\P}(D)$ be an invertible sheaf on $\P$ corresponding to an ample
divisor $D$, and $\Delta$ be the support polytope for $\O_{\P}(D)$. Then
\begin{itemize}
\item[{\rm 1)}] 
  \[
    h^q(\P,\Omega_{\P}^p) =
    \left\{
    \begin{array}{ll}
       \displaystyle
       \sum_{s = p}^n(-1)^{s + p}\binom{s}{p} f_s,  &q = p,\\
       \\
       0,                                           &q\ne p.
    \end{array}
    \right.
\]
\item[{\rm 2)}]
  \[
    h^q(\P,\Omega_{\P}^p(D)) =
    \left\{
    \begin{array}{ll}
       \displaystyle
       \sum_{s = p}^n\binom{s}{p} 
       \sum_{\Gamma\in\F_s(\Delta)} \l^*(\Gamma),  &q = 0,\\
       \\
       0,                                          &q > 0.
    \end{array}
    \right.
\]
\end{itemize}

\label{Bott.formula.II}
\end{theo}
\begin{proof} 
The formula
\[
  h^p(\P,\Omega_{\P}^p) =
  (-1)^p\sum_{\Gamma\subseteq \Delta}(-1)^{\dim\Gamma}\binom{\dim\Gamma}{p} =
  \sum_{s = p}^n(-1)^{s + p}\binom{s}{p} f_s
\]
is the assertion of Corollary 2.5 of \cite{Danilov-Khovanskii}.
We prove this statement independently.
Making use of the exact sequence (\ref{restr.hom}), we get
\[
  \chi(\P,\Omega_\P^p) = 
  \sum_{\Gamma\subseteq\Delta}
  \chi(\P_{\Gamma},\Omega_{(\P_{\Gamma},\,D_{\Gamma})}^p).
\]
Now the requested formula follows from Proposition \ref{Remark.in.DK}
and the vanishing of $h^q(\P,\Omega_{\P}^p)$ for $p \ne q$.

Let us prove assertion $2)$.
The short exact sequence (\ref{restr.hom}) remains exact after 
tensoring by $\L = \O_{\P}(D)$.
As in the case above, we have
\[
  \chi(\P,\Omega_\P^p\otimes\L) = 
  \sum_{\Gamma\subseteq\Delta}
  \chi(\P_{\Gamma},\Omega_{(\P_{\Gamma},\,D_{\Gamma})}^p\otimes\L),
\]
and from Proposition \ref{Remark.in.DK} and the Bott vanishing
theorem,
\[
  \chi(\P,\Omega_\P^p\otimes\L) = 
  h^0(\P,\Omega_\P^p\otimes\L)  = 
  \sum_{\Gamma\subseteq\Delta}
  \binom{\dim\Gamma}{p}\l^*(\Gamma) =
  \sum_{s = p}^n\binom{s}{p} 
  \sum_{\Gamma\in\F_s(\Delta)} \l^*(\Gamma).
\]

We  prove the second assertion of the theorem by using the description 
of the space of global sections of the sheaf $\Omega_{\P}^p\otimes \L$.
We have a decomposition into a direct sum of $M$-homogeneous components

\[
  H^0(\P,\Omega_{\P}^p\otimes \L) = \bigoplus_{m\in \Delta\cap M}
  \bigwedge\nolimits^p V_\Delta (m),
\]
where $V_\Delta(m)$ is the $\C$-subspace in $M_\C:=M\otimes\C$ generated by the
smallest face of $\Delta = \Delta(\L)$ containing $m$.
Since the support polytope $\Delta$ admits the natural partition
$\Delta = \coprod_{\Gamma} Int(\Gamma)$, where $Int(\Gamma)$ is 
the relative interior of the face $\Gamma$, we have
\[
  h^0(\P,\Omega_{\P}^p\otimes \L)=
  \sum_{\Gamma\subseteq \Delta}
  \sum_{m\in Int(\Gamma)\cap M}
  \binom{\dim V_{\Delta}(m)}{p}=
  \sum_{\Gamma\subseteq \Delta}
  \binom{\dim\Gamma}{p} \l^*(\Gamma).
\]
Thus we have given a different proof for the second formula.
\end{proof}
\begin{coro} One has the following equality
  \begin{eqnarray*}
    \displaystyle
    \sum_{p = 0}^n h^0(\P,\Omega_{\P}^p(D)) y^p = 
    \displaystyle
    \sum_{s = 0}^n \Delta_{(s)}
    (y + 1)^k,\quad \mbox{\it if} \; D \; \mbox{\it is ample},
  \end{eqnarray*}
where $\Delta_{(s)} := \displaystyle\sum_{\Gamma\in\F_s(\Delta)}\l^*(\Gamma)$.
\label{Bott.Corollary.II}
\end{coro}

\section{Relation with the original Bott formula and some \\ 
combinatorial identities}
In this section we compare our results of theorems \ref{Bott.formula.I}
and \ref{Bott.formula.II} with the original Bott formula of 
Theorem \ref{Original.Bott}. 
Recall that $\O_{\P^n}(k)$ is ample if and only if $k>0$, and the support 
polytope for $\O_{\P^n}(k)$ is the simplex
\begin{eqnarray*}
  \Delta = \Delta(k) =
  \{m\in M_\R:\, m_1\ge 0,\ldots,m_n\ge 0,\quad m_1+\ldots+m_n\le k\}.
\end{eqnarray*}
It is easy to see that
\[
  f_{n - j} = \binom{n+1}{j},\quad
  \Delta^{(j)} = \binom{n+k-j}{k}\binom{n+1}{j},
\]
and
\[
  f_s = \binom{n+1}{n-s},\quad
  \Delta_{(s)} =
  \binom{k-1}{s}\binom{n+1}{n-s},
\]
for all $0\le j\le n$ and $0\le s\le n$.
Comparison of the formulas $1)$ of Theorems \ref{Bott.formula.I} and
\ref{Original.Bott} corresponding to the case $q = p$ gives the identity
$$
  \sum_{j=0}^p(-1)^{p-j}\binom{n-j}{p-j}\binom{n+1}{j} = 1.
  \eqno{(a')}
$$
Now we compare the formulas $2)$ in these theorems corresponding to the 
case $q = 0$:
$$
  \sum_{j=0}^p(-1)^j\binom{n-j}{p-j}\binom{n+k-j}{k}\binom{n+1}{j}=
  \binom{n+k-p}{k}\binom{k-1}{p}.
  \eqno{(b')}
$$
Analoguous comparisons of Theorems \ref{Bott.formula.II} and
\ref{Original.Bott} show that
$$
  \sum_{s=p}^n(-1)^{p+s}\binom{s}{p}\binom{n+1}{n-s}=1,
  \eqno{(a'')}
$$
and
$$
  \sum_{s=p}^n\binom{s}{p}\binom{k-1}{s}\binom{n+1}{n-s}=
  \binom{n+k-p}{k}\binom{k-1}{p}.
  \eqno{(b'')}
$$
These identities can be proved directly.

Let us verify the simplest identity $(a')$.
To prove $(a')$ it is sufficient to check the functional equation
\[
  \sum_{j=0}^n\binom{n+1}{j}t^j(1-t)^{n-j}=\sum_{k=0}^n t^k,
\]
and then to compare the coefficients by the monomial $t^p$.
This equation follows from Corollary \ref{Bott.Corollary.I},
but we can deduce it from the obvious identity
\[
  (t+1-t)^{n+1}=1,
\]
or
\[
  (1-t)\sum_{j=0}^n\binom{n+1}{j}t^j(1-t)^{n-j}+t^{n+1}=1.
\]
Now we get the required by $\displaystyle (1 - t^{n+1})(1 - t)^{-1}=\sum_{k = 0}^n t^k$.

The identity $(a'')$ is equivalent to $(a')$.
The remaining identities $(b')$ and $(b'')$ have been proved in 
\cite{Materov-Yuzhakov}.
Our proof is based on the method of integral representations for 
combinatorial sums in the spirit of the work \cite{Egorychev-Yuzhakov}.

\section{The Hilbert-Ehrhart polynomials and 
the Serre duality}


Let $\Delta\subset\R^n$ be a non-empty integer polytope. 
For a positive integer $k$, let $k\cdot\Delta := \{kx:\,x\in\Delta\}$ denote
the dilated polytope. It was proved by Ehrhart \cite{Ehrhart}, and 
in somewhat stronger form by Macdonald \cite{Macdonald},
that there is a polynomial $L(k)$ with the properties 
\begin{itemize}
\item[(i)] for any integer $k > 0$, one has
\[
  L(k) = \l(k\cdot\Delta);
\]
\item[(ii)] for any integer $k > 0$, one has the following 
{\it reciprocity law}
\[
  L(-k) = \l^*(k\cdot\Delta).
\]
\end{itemize}
The polynomial $L(k)$ is called the {\it Ehrhart polynomial} for $\Delta$.
In this section we study a generalization of the Ehrhart polynomial 
using the techniques of algebraic geometry (cf., \cite{Danilov, Oda.book}). 
The basic idea is to compare the Ehrhart polynomial with a certain 
Hilbert polynomial.
For example, the Serre duality provides a generalization of the
reciprocity law.

\begin{prop}[\cite{Snapper, Kleiman}] 
Let $\L$ be an invertible sheaf 
on a complete variety $\P$ and $\F$ be a coherent sheaf of $\O_{\P}$-modules 
on $\P$. Then the Euler-Poincar\'e characteristic 
$\chi(\P,\F\otimes \L^{\otimes k})$ is a polynomial in $k$  
of total degree $\le n$ which assumes integer values whenever $k$ is
integer.
\end{prop}
\begin{dfn} 
Let $\L=\O_{\P}(D)$ be an invertible sheaf corresponding to
the ample Cartier divisor $D$ of a complete simplicial toric variety $\P$.
We call the polynomial
\[
  L_p(k) := \chi(\P,\Omega_{\P}^p\otimes \L^{\otimes k})=
  \chi(\P,\Omega_{\P}^p(k D))
\]
the {\it Hilbert polynomial} of the sheaf $\Omega_{\P}^p$ with respect
to $\L$, or {\it $p$-th Hilbert polynomial}.
\end{dfn}

It follows from the Bott formula and Bott vanishing theorem 
for toric varieties, that
\begin{eqnarray*}
  L_p(k) = h^0(\P,\Omega_{\P}^p(k D)) 
  &=& 
  \sum_{j = 0}^p (-1)^j \binom{n-j}{p-j} 
  \sum_{F \in \F_{n - j}(\Delta)} \l(k\cdot F) \\
  &=& 
  \sum_{s = p}^n \binom{s}{p} 
  \sum_{\Gamma \in \F_{s}(\Delta)} \l^*(k\cdot \Gamma)
\end{eqnarray*}
for $k > 0$. For example, for $p = 0$ the polynomial $L_0(k)$ coincides
with the usual Ehrhart polynomial of the polytope $\Delta = \Delta(\L)$,
and with $L_n(k) = \l^*(k\cdot \Delta)$ whenever $k > 0$.
Hence, we can consider $L_p(k)$ as a generalization of the Ehrhart
polynomial and also call it the {\it $p$-th Ehrhart polynomial} for 
$\Delta$, or simply the {\it $p$-th Hilber-Ehrhart polynomial}. 

Note that the reciprocity law can be written as
\[
  L_0(-k) = (-1)^n L_n(k).
\]
Here we prove a more general result.
\begin{theo} The polynomial $L_p(k)$ satisfies the duality property
\[
  L_p(-k)=(-1)^n L_{n-p}(k)
\]
for any positive integer $k$ and $0 \le p \le n$.
\label{gen.rec.law}
\end{theo}
\begin{proof}
We need a special form of the Serre-Grothendieck duality 
(see \S 3.3 of \cite{Oda.book})
\[
  H^q(\P,\Omega_{\P}^p\otimes\F)^* \simeq
  H^{n-q}(\P,\Omega_{\P}^{n - p}\otimes\F^*),
  \quad 0 \le q \le n
\]
for any locally free sheaf of the $\O_{\P}$-modules $\F$ with the dual
$\F^*$ on a complete simplicial toric variety $\P$.
Take $\F = \O_{\P}(kD)$, $k > 0$. 
From the Serre duality we have the isomorphism
\[
  H^0(\P,\Omega_{\P}^p(k D))^* \simeq
  H^n(\P,\Omega_{\P}^{n-p}(- k D)).
\]
Hence, using the vanishing of $h^q(\P,\Omega_{\P}^p(k D))$ for $q > 0$,
we get the equality
\[ 
  \chi(\P,\Omega_{\P}^p(k D)) = 
  (-1)^n\chi(\P,\Omega_{\P}^{n-p}(-k D)),
\]
which is equivalent to 
\[ 
  \chi(\P,\Omega_{\P}^p(-k D)) = 
  (-1)^n\chi(\P,\Omega_{\P}^{n-p}(k D)),
\]
and this completes the proof.
\end{proof}

Now let $\P$ be nonsingular.
Recall, that any variety $\P$ has a Todd homology class $td_{\P}$ of $\P$ 
in $A_*(\P)_{\Q}$, see \cite{Fulton}. 
Since $\P$ is nonsingular, $td_{\P}=Td_{\P}\cap [\P]$,
where $Td_{\P}$ is the Todd cohomogy class in $A^*(\P)_{\Q}$ 
and $[\P]$ is the fundamental class of $\P$. 

We state without proof some facts on intersection theory on toric 
varieties, where part ${\rm (a)}$ slightly generalizes the results in
\cite{Danilov, Fulton, Oda.book}, and ${\rm (b)}$ obviously follows from the
Hirzebruch-Riemann-Roch theorem \cite{Hirzebruch}.
\begin{prop} Let $L_p(k)$ be the Hilbert-Ehrhart polynomial.

{\rm (a)} The coefficients of the polynomial 
$L_p(k) = \sum_{i = 0}^n a_{pi} k^i$
are the intersection numbers
\[
  a_{pi} = 
  \displaystyle\frac{1}{i!}\int_{\P}\,D^i\cdot ch(\Omega_{\P}^p)\cdot Td_{\P},
\]
where $\int_{\P}$ denotes the degree homomorphism 
$\int_{\P}:A^0(\P)_\Q\,\rightarrow \Q$ and $ch(\Omega_{\P}^p)$ is the Chern 
character of the sheaf $\Omega_{\P}^p$. For example, the leading coefficient is
\[
  \displaystyle
  a_{pn} = \binom{n}{p}Vol(\Delta),
\]
where $Vol(\Delta)$ is the normalized volume of the polytope 
$\Delta\subset M_{\R}$ such that the volume of the unit cube determined by
the basis of the lattice $M$ is $1$.

{\rm (b)} The generating function for the polynomial $L_p(k)$ is equal to 
\[
  \displaystyle
  L(y;k) := \sum_{p = 0}^n L_p(k) y^p =
  \int_{\P} e^{kD(y + 1)}\prod_{j = 1}^d Q(y;D_j),
\]
where 
\[
  Q(y;x) = \frac{x(y + 1)}{e^{x(y + 1)} - 1} + x,
\]
and $D_j$ are the torus-invariant divisors corresponding to the 
one-dimensional generators of the fan $\Si$.
\end{prop}

\section{Combinatorics of simple polytopes}

Comparison of two versions of the Bott formula gives some 
elegant corollaries in combinatorics of simple polytopes.

If we compare the first parts of Theorems \ref{Bott.formula.I} and
\ref{Bott.formula.II}, we get the well-known Dehn-Sommerville equations,
conjectured by Dehn \cite{Dehn} and proved by Sommerville \cite{ Sommerville}.
(See also \cite{Broendsted}.) In algebro-geometrical context these equalities
were proved in \cite{Stanley} for any simple rational polytope, using the 
Poincar\'e duality on toric varities, and by Oda \cite{Oda.book}, using 
the Serre-Grothendieck duality theorem.

\begin{theo} For any simple lattice $n$-polytope $\Delta$ 
and $p = 0,1,\ldots,n$ one has the identities
\[
 \sum_{j = 0}^p (-1)^j \binom{n-j}{p-j} f_{n-j}
  =
 \sum_{s = p}^n (-1)^{s}\binom{s}{p} f_s.
\]
\end{theo}

Now compare the second parts of Theorems \ref{Bott.formula.I} and
\ref{Bott.formula.II}. We get non-trivial relations between the integer 
points in faces of the simple integer polytope $\Delta$, which we prove 
here in a purely combinatorial way.
\begin{theo} For any simple lattice $n$-polytope $\Delta$ 
and $p = 0,1,\ldots,n$ one has the identities
\label{comb_dual}
\[
 \sum_{j = 0}^p (-1)^j \binom{n-j}{p-j} 
 \sum_{F \in \F_{n - j}(\Delta)} \l(F) =
 \sum_{s = p}^n \binom{s}{p} 
 \sum_{\Gamma \in \F_{s}(\Delta)} \l^*(\Gamma).
\]
\end{theo}
\begin{proof} 
We have the chain of equalities
\begin{eqnarray*}
 \sum_{s = p}^n \binom{s}{p} 
 \sum_{\Gamma \in \F_{s}(\Delta)} \l^*(\Gamma) 
 &=&
 \sum_{\Gamma \subseteq \Delta} \binom{\dim\Gamma}{p} \l^*(\Gamma) 
 \\
 &=&
 \sum_{\Gamma \subseteq \Delta} \binom{\dim\Gamma}{p} 
 \sum_{F \subseteq \Gamma} (-1)^{\dim\Gamma - \dim F}\l(F)
 \\
 &=& 
 \sum_{F \subseteq \Delta} (-1)^{\dim F} \l(F) 
 \sum_{\Gamma \supseteq F} (-1)^{dim \Gamma} \binom{\dim\Gamma}{p}
 \\
 &=&
 \sum_{F \subseteq \Delta} (-1)^{n + \dim F} \binom{\dim F}{n - p} \l(F)
 \\
 &=&
 \sum_{j = n - p}^n (-1)^{n + j} \binom{j}{n-p}
 \sum_{F \in \F_j(\Delta)} \l(F)
 \\
 &=&
 \sum_{j = 0}^p (-1)^j \binom{n-j}{p-j} 
 \sum_{F \in \F_{n - j}(\Delta)} \l(F),
\end{eqnarray*}
where we have used the identity
\[ 
 \sum_{\Gamma \supseteq F} (-1)^{dim \Gamma} \binom{\dim\Gamma}{p}
 = 
 (-1)^{n} \binom{\dim F}{n - p},
\] 
which is proved in the Appendix.
\end{proof}
\begin{rem} 
There are two obvious cases of the identities of Theorem \ref{comb_dual}.
For $p = 0$, we have
\[
  \l(\Delta) = \sum_{\Gamma\subseteq\Delta} \l^*(\Gamma).
\]
For $p = n$, we get the inclusion-exclusion formula
\[
  \l^*(\Delta) = \sum_{j = 0}^n (-1)^{n - j} \sum_{F\in\F_{j}(\Delta)}\l(F).
\]
\end{rem}

We are in a position to give a purely combinatorial proof of the generalized 
reciprocity law of Theorem \ref{gen.rec.law}.
\begin{coro} 
For any simple lattice $n$-polytope $\Delta$ and any positive integer 
$k$, one has
\[
 L_p(-k) = (-1)^n L_{n-p}(k),
\]
for any $p = 0,\ldots,n$.
\end{coro}
\begin{proof} Using the equality of Theorem \ref{comb_dual}, we get
\begin{eqnarray*}
 L_p(-k) 
 &=& 
 \sum_{j = 0}^p (-1)^j \binom{n-j}{p-j} 
 \sum_{F \in \F_{n - j}(\Delta)} (-1)^{n - j} \l^*(k\cdot F) \\
 &=& 
 (-1)^n  
 \sum_{s = p}^n \binom{s}{n - p} 
 \sum_{F \in \F_{s}(\Delta)} \l^*(k\cdot F) \\
 &=&
 (-1)^n L_{n-p}(k),
\end{eqnarray*}
as required.
\end{proof}

\begin{exam} Let $\Delta$ be a convex polytope in $\R^3$ with vertices at
$(0,\,0,\,0)$, $(1,\,0,\,0)$, $(0,\,1,\,0)$ and $(1,\,1,\,m)$, where $m$
is a positive integer. Then
\begin{eqnarray*}
  \begin{array}{l}
    \displaystyle
    L_0(k)=\frac{m}{6}k^3 + k^2 + \frac{12-m}{6}k + 1,\\
    \\
    \displaystyle
    L_1(k)=\frac{m}{2}k^3 + k^2 - \frac{m}{2}k - 1,\\
    \\
    \displaystyle
    L_2(k)=\frac{m}{2}k^3 - k^2 - \frac{m}{2}k + 1,\\
    \\
    \displaystyle
    L_3(k)=\frac{m}{6}k^3 - k^2 + \frac{12-m}{6}k - 1.
  \end{array}
\end{eqnarray*}
\end{exam}


\section{Weighted components of sheaves with 
logarithmic poles}

Recall that if $\L$ is an ample invertible sheaf on a complete simplicial
toric variety $\P$, then
\[
  H^q(\P,\Omega_\P^p\otimes \L) = 0
\]
for all $q>0$ and $p\ge 0$ (see \cite{Danilov}). 
This statement is a generalization of the theorems of Bott and Steenbrink 
for projective spaces and weighted projective spaces respectively.
Batyrev and Cox proved a more general vanishing theorem. 

\begin{dfn}[\cite{Batyrev-Cox}]
Denote by $\Omega_\P^p({\rm log} (-K))$ the sheaves of differential $p$-forms 
with logarithmic poles along the anticanonical divisor $-K$ on $\P$. 
Let $\mathcal{W}$ be the {\it weight filtration}
\[
  \mathcal{W}:\, 0\subset
  W_0\Omega_{\P}^p({\rm log} (-K))\subset
  W_1\Omega_{\P}^p({\rm log} (-K))\subset\cdots\subset
  W_p\Omega_{\P}^p({\rm log} (-K))=\Omega_{\P}^p({\rm log} (-K))
\]
on $\Omega_{\P}^p({\rm log} (-K))$ defined by
$W_k\Omega_{\P}^p({\rm log} (-K)):=
\Omega_{\P}^{p-k}\wedge\Omega_{\P}^k({\rm log} (-K))$.
\end{dfn}
\begin{theo}[\cite{Batyrev-Cox}, Theorem 7.2] 
Let $\L$ be an ample invertible sheaf on a complete simplicial
toric variety $\P$. Then for any $p\ge 0$, $k\ge 0$, and $q>0$, one has
\[
  H^q(\P,W_k\Omega_{\P}^p({\rm log} (-K))\otimes \L)=0.
\]\label{B.S.D}
\end{theo}
Here we give some additional information about the global sections of
$W_k\Omega_{\P}^p({\rm log} (-K))\otimes \L$. We need the following 
result of V.I.~Danilov.
\begin{theo}[\cite{Danilov}, \S 15.7] 
For any integer $0\le k\le p$ one has the short exact sequence
   \begin{equation}
     0\rightarrow W_{k-1}\Omega_{\P}^p({\rm log} (-K))
       \rightarrow W_k   \Omega_{\P}^p({\rm log} (-K))
       \stackrel{\rm Res}{\rightarrow}
       \bigoplus_{\s\in\Si(k)}\Omega_{V(\s)}^{p-k}
       \rightarrow 0,
   \label{weighted.sequence}
   \end{equation}
where $V(\s)$ is the Zariski closure of the torus-invariant orbit of $\P$
corresponding to $\s\in\Si(k)$ and {\rm Res} is the Poincar\'e residue map.
\end{theo}
It is easy to see that by tensoring by $\L$ the short exact sequence
(\ref{weighted.sequence}) and applying the vanishing Theorem 
\ref{B.S.D} we obtain
  \begin{eqnarray*}
    0\rightarrow H^0(\P,W_{k-1}\Omega_{\P}^p({\rm log} (-K))\otimes \L)
     \rightarrow H^0(\P,W_k    \Omega_{\P}^p({\rm log} (-K))\otimes \L)
     \rightarrow\\
     \rightarrow
     \bigoplus_{\s\in\Si(k)}H^0(\P,\Omega_{V(\s)}^{p-k}\otimes \L)
     \rightarrow 0.
  \end{eqnarray*}
Thus, using an induction, we get
\[
  h^0(\P,W_k \Omega_{\P}^p({\rm log} (-K))\otimes \L)=
  \sum_{s=0}^k\sum_{\s\in\Si(s)}
  h^0(\P,\Omega_{V(\s)}^{p-s}\otimes \L).
\]
The generalized Bott formula from Theorem \ref{Bott.formula.I} and Proposition
\ref{support.polytope} implies
\[
  h^0(\P,\Omega_{V(\s)}^{p-s}\otimes \L)=
  \sum_{j=0}^{p-s}(-1)^j\binom{n-s-j}{p-s-j}
  \sum_{\begin{subarray}{l}
           \t\in\Si(j) \\
           \t\succ\s
        \end{subarray}}
  \l(\Delta_\t).
\]
We have proved the following statement.
\begin{prop} Let $\P$ be a complete simplicial toric variety and
$\L$ be an ample invertible sheaf on $\P$ with the support polytope
$\Delta$. Then one has the equality
\[
  h^0(\P,W_k \Omega_{\P}^p({\rm log} (-K))\otimes \L)=
  \displaystyle
  \sum_{s=0}^k
  \sum_{\s\in\Si(s)}
  \sum_{j=0}^{p-s}(-1)^j\binom{n-s-j}{p-s-j}
  \sum_{\begin{subarray}{l}
           \t\in\Si(j) \\
           \t\succ\s
        \end{subarray}}
  \l(\Delta_\t).
\]
\end{prop}

\section{Cohomology of quasi-smooth hypersurfaces
\label{section.cohomologies}
}

In this section we reprove the ``combinatorial part'' of the well-known result
of Danilov and Khovanski\^i on the Hodge numbers of quasi-smooth 
hypersurfaces in toric varieties. 
Our proof is simple and uses the Bott formula for toric varieties.

Any complete simplicial toric variety $\P = \P(\Si)$ can be constructed as a 
geometric quotient as follows (see \cite{Cox}). 
Suppose that $S(\Si) := \C[z_1,\ldots,z_d]$ is the 
polynomial ring over $\C$ with the variables $z_i$ corresponding to the 
integer generators $v_i$ of the one-dimensional cones of $\Si$. 
For every $\s\in\Si$ let
$\hat{z}_\s := \prod_{v_i\not\in\s} z_i$, and let
\[
  Z(\Si) := \bigcap_{\s\in\Si}\{z\in\C^d:\,\hat{z}_\s=0\}.
\]
The toric variety $\P$ is a geometric quotient of the Zariski open set
$U(\Si) := \C^d\backslash Z(\Si)$ by the algebraic group
$G(\Si) := {\rm Hom}_\Z(A_{n-1}(\P),\,\C^*)$, where $A_{n-1}(\P)$ is the Chow
group of $\P$. The action of $G(\Si)$ on $\C^d$ induces the grading on 
$S(\Si)$.
Moreover, if $\L$ is an invertible sheaf on $\P$, then for
$\alpha = [\L]\in A_{n-1}(\P)$ one has the isomorphism 
$H^0(\P,\L)\simeq S_\alpha$.
A polynomial $f$ in the graded part $S_\alpha$ of $S(\Si)$ corresponding
to the class $\alpha\in A_{n-1}(\P)$ is said to be $G$-{\it homogeneous}
of degree $\alpha$. Thus, the global sections of $\L$ determine the Cartier
divisor $D$ on $\P$.
\begin{dfn}[\cite{Batyrev-Cox}]
The hypersurface $D\subset \P$ defined by the $G$-homogeneous polynomial $f$
in $S(\Si)$ is said to be {\it quasi-smooth} if ${\bf V}(f)\cap U(\Si)$
is either empty or a nonsingular subvariety of codimension one in $U(\Si)$.
\end{dfn}

\begin{dfn} 
Let $X$ be an algebraic variety over $\C$.
Denote by $h^{p,\,q}(H^k(X,\C))$ the dimension of the $(p,\,q)$-component 
of the $k$-th cohomology group.
Let us introduce the numbers
\begin{eqnarray*}
  &&e^{p,\,q}(X) := \sum_k (-1)^k h^{p,\,q}(H^k(X,\C)),\\
  &&e^p (X)      := \sum_q e^{p,\,q}(X).
\end{eqnarray*}
For a complete nonsingular variety $X$, we have 
\[
  e^{p,\,q}(X) = (-1)^{p + q} h^{p,\,q}(X),
\]
where $h^{p,\,q}(X)$ are the (ordinary) Hodge numbers of $X$.
\end{dfn}
\begin{rem} 
Danilov and Khovanski\^i use cohomologies with compact supports
instead of usual cohomologies, but in our situation there exists 
the Poincar\'e duality
\begin{eqnarray}
  h^{p,\,q}(H^k(X,\C)) =  h^{r-p,\,r-q}(H_c^{2r-k}(X,\C)),
\label{Poinc}
\end{eqnarray}
where $r = \dim_{\C}X$, which relates the two treatments.
\end{rem}

Let $D\subset \P$ be a nondegenerate divisor on a complete simplicial
$n$-dimensional toric variety. Then the natural map
\[
  H^i(\P)\rightarrow H^i(D)
\]
is an isomorphism for $i < n - 1$ and is injective for $i = n - 1$.
\begin{dfn} 
Define the {\it primitive cohomology} $H_0^{n - 1}(D)$ by
the exact sequence
\[
  0 \rightarrow H^{n - 1}(\P)
    \rightarrow H^{n - 1}(D)
    \rightarrow H_0^{n - 1}(D)
    \rightarrow 0.
\]
\end{dfn}
\begin{rem} 
As A.~Mavlyutov pointed out to us (see 
\cite{Mavlyutov2}, Remark 5.5), the primitive cohomologies $H_0^{n - 1}(D)$
coincide with the residue part $H_{res}^{n - 1}(D)$ of cohomology defined 
as the image of the residue map 
$Res:\,H^n(\P\backslash D)\rightarrow H^{n - 1}(D)$.
\end{rem}

Recall, that since $\P$ is simplicial and $D$ is nondegenerate, 
$H_0^{n - 1}(D)$ has a pure Hodge structure (cf., \cite{Batyrev-Cox}).
\begin{theo}[\cite{Danilov-Khovanskii}] 
Let $\Delta$ be a support polytope, corresponding
to an ample nondegenerate hypersurface $D$. Then,
\[
  h_0^{p,\,n - 1 - p}(D) = 
  (-1)^n \sum_{\Gamma\subseteq\Delta} 
  \varphi_{\dim\Gamma - p}(\Gamma),
\]
where
\[
  \varphi_i(\Gamma) := 
  \sum_{k = 1}^i (-1)^{i - k} \binom{\dim\Gamma + 1}{i - k} \l^*(k\cdot\Gamma).
\]
\label{D-Kh}
\end{theo}

Denote by $\Omega_{\P}^p(\logD)$ the sheaf of $p$-differential forms 
on $\P$ with logarithmic poles along $D$ (see \S 15 of \cite{Danilov}).
First, we compute the Euler-Poincar\'e characteristic
of the sheaf $\Omega_{\P}^p(\logD)$.
\begin{prop}[\cite{Batyrev-Cox}, Proposition 10.1] 
If $D$ is a quasi-smooth hypersurface of a complete toric variety
$\P$ defined by zeros of a global section of the ample invertible sheaf 
$\O_{\P}(D)$ on $\P$, then there is an exact sequence
  \begin{eqnarray*}
    0\rightarrow \Omega_{\P}^p(\logD)
    \rightarrow
    \Omega_{\P}^p(D) \stackrel{d}{\rightarrow}
    \Omega_{\P}^{p+1}(2D)/\Omega_{\P}^{p+1}(D)
    \stackrel{d}{\rightarrow}
    \ldots\\
    \ldots\stackrel{d}{\rightarrow}
    \Omega_{\P}^n((n-p+1) D)/\Omega_{\P}^n((n-p) D)
    \rightarrow 0.
  \end{eqnarray*}
\label{Long.exact.sequence}
\end{prop}
\begin{lem} Let $\Delta$ be the support polytope for $\O_{\P}(D)$
corresponding to the ample divisor $D$. Then the Euler-Poincar\'e 
characteristics of the sheaf $\Omega_{\P}^p(\logD)$ is equal to
\[
  \chi(\P,\Omega_{\P}^p(\logD))=
  \sum_{k = 1}^{n - p}(-1)^{k + 1}
  \sum_{\Gamma\subseteq \Delta}
  \binom{\dim\Gamma + 1}{p + k}\l^*(k\cdot\Gamma).
\]
\end{lem}
\begin{proof}
Indeed, from the long exact sequence of  
Proposition \ref{Long.exact.sequence}, we have
\[
  \chi(\P,\Omega_{\P}^p(\logD))=
  \sum_{k = 1}^{n - p}
  (-1)^{k + 1}
  [\chi(\P,\Omega_{\P}^{p + k - 1}(kD)) +
   \chi(\P,\Omega_{\P}^{p + k}(kD))].
\]
Since $D$ is ample, 
\[
  \chi(\P,\Omega_{\P}^p(\logD))=
  \sum_{k = 1}^{n - p}
  (-1)^{k + 1}
  [h^0(\P,\Omega_{\P}^{p + k - 1}(kD)) +
   h^0(\P,\Omega_{\P}^{p + k}(kD))],
\]
applying the generalized Bott formula given in Theorem 
\ref{Bott.formula.II}, we prove the statement.
\end{proof}

We start with the relation 
\[
  e^p(D) = e^{p + 1}(\P) - e^{p + 1}(\P\backslash D),
\]
following from the Gysin exact sequence
\[
  \ldots\, \rightarrow \, H^{k-2}(D)          \, 
           \rightarrow \, H^k(\P)             \, 
           \rightarrow \, H^k(\P\backslash D) \, 
           \rightarrow \, H^{k - 1}(D)        \, 
           \rightarrow \, 
  \ldots
\]
of hypercohomologies of the exact sequence of complexes
\[
  0\, \rightarrow \, \Omega_{\P}^{\cdot}        \,
      \rightarrow \, \Omega_{\P}^{\cdot}(\logD) \,
      \rightarrow \, \Omega_{D}^{\cdot - 1}     \,
      \rightarrow \, 0.
\]
There are two spectral sequences
\[
  'E_1^{p,\,q}  = H^q (\P,\Omega_{\P}^p) 
  \Rightarrow H^{p + q}(\P,\C),
\]
and
\[
  ''E_1^{p,\,q} = H^q (\P,\Omega_{\P}^p(\logD))
  \Rightarrow H^{p + q}(\P\backslash D,\C),
\]
both degenerating at the first term and converging to the Hodge 
filtration $\F^\cdot$ on $H^*(\P,\C)$ and $H^*(\P\backslash D,\C)$ 
respectively \cite{Danilov}. Therefore,
we have the equalities
\[
  e^{p + 1}(\P) = (-1)^{p + 1} \chi(\P,\Omega_{\P}^{p + 1}),
  \quad
  e^{p + 1}(\P\backslash D) = 
  (-1)^{p + 1} \chi(\P,\Omega_{\P}^{p + 1}(\logD)).
\]
Using the  Bott formula, one obtains:
\begin{prop} One has the equality
\[
  e^p(D) = 
  (-1)^{p + 1} \sum_{k \ge 0} 
  (-1)^k \binom{k}{p + 1} f_k -
  \sum_{\Gamma\subseteq \Delta}
  (-1)^{\dim \Gamma}
  \varphi_{\dim\Gamma - p} (\Gamma).
\]
\end{prop}

By applying the Poincar\'e duality (\ref{Poinc}), we restore all the 
Hodge numbers $h_0^{p,\,q}(D)$ and get the statement of Theorem \ref{D-Kh}.

\section{The Bott formula on $\P(\E)$}

Let $\L_0,\ldots,\L_s$ be ample invertible sheaves on a complete 
$n$-dimensional toric variety $\P$. Denote $Y = \P(\E)$ 
the projective space bundle associated with the sheaf 
$\E := \L_0\oplus\cdots\oplus \L_s$, with the invertible sheaf $\O_Y(1)$ 
and the projection $\pi:\,Y\rightarrow \P$. The sheaf $\O_Y(1)$ is ample since
the $\L_j$ are ample (see \S1 of \cite{Hartshorn}).
\begin{dfn}
Let $\Omega^1_{Y/\P}$ be the sheaf of relative differentials arising from 
the short exact sequence
(see \cite{Manin})
\begin{eqnarray}
\label{Omega_Y}
 0 \rightarrow \Omega^1_{Y/\P} 
   \rightarrow \pi^*(\E)\otimes \O_Y(-1)
   \rightarrow \O_Y
   \rightarrow 0.
\end{eqnarray}
\end{dfn}
The purpose of this section is to give a combinatorial description of 
the global sections of the sheaf of relative $p$-th differentials
$\Omega^p_{Y/\P}(k) := 
\bigwedge^p \Omega^1_{Y/\P}\otimes \O_Y(k)$.

Recall the construction of $Y$ as a toric variety (see \cite{Oda.book}).
Suppose that the support polytope associated with the sheaf $\L_j$ has the form
\[
 \Delta_j = \{ m\in M_{\R}: \<m,v_i\>\ge -a_{ij},\, i = 1,\ldots,d\}.
\]
Let $N'\simeq \Z^s$ be a lattice with $\Z$-basis $\{e_1,\ldots,e_s\}$
and $\widetilde{N} := N \oplus N'$. The $1$-dimensional cones of the fan
$\widetilde{\Sigma}$ corresponding to the toric variety $Y$ have the 
generators
\begin{eqnarray*}
 \widetilde{v}_i &=& v_i + \sum_{j = 1}^s (a_{ij} - a_{i0}) e_j, 
                      \quad i = 1,\ldots,d,\\
 \widetilde{e}_0 &=& -e_1 - \ldots -e_s, \\
 \widetilde{e}_j &=& e_j, \quad j = 1,\ldots,s.
\end{eqnarray*}
Denote by $\widetilde{\sigma}$ the image of each $\sigma\in\Sigma$ 
under the map $N_{\R}\rightarrow \widetilde{N}_{\R}$ given by 
$v_i\mapsto \widetilde{v}_i$ and let $\sigma'$ be the cones generated
by ${e}_0,\ldots,{e}_{i-1},{e}_{i+1},\ldots,{e}_s$. Then the fan 
$\widetilde{\Sigma}$ is generated by the cones $\widetilde{\sigma} + \sigma'$
and their faces.
Let $\eta_1,\ldots,\eta_s$ be the generators of the lattice $M'$
dual to $N'$. Denote the torus invariant divisors on $Y$ corresponding
to $\widetilde{v}_i$ and $\widetilde{e}_j$ by $\widetilde{D}_i$ and 
$\widetilde{D}_j'$ respectively. Note that 
$\pi^*(\L_j) = \O_Y(\sum_{i = 1}^d a_{ij} \widetilde{D}_i)$.
\begin{lem} The polytope
\label{support_O_Y}
\[
 \nabla = 
 \{
 \lambda_1 \eta_1 + \ldots + \lambda_s \eta_s +
 \lambda_0 m_0 + \ldots + \lambda_s m_s :\, \lambda_j\ge -k_j,\,
 \sum_{j = 0}^s \lambda_j = k,\, m_j \in \Delta_j
 \}
\]
in $\widetilde{M}_{\R} := M_{\R}'\oplus M_\R$ is the support polytope 
associated with the sheaf 
$\O_Y(k)\otimes
\O_Y(k_0\widetilde{D}_0' + \ldots + k_s \widetilde{D}_s')$.
\end{lem}
\begin{proof} It follows from the Lemma 2.1 in \cite{Mavlyutov} that
\[
  \O_Y(1) \simeq \O_Y(\widetilde{D}_0')\otimes\pi^*(\L_0) =
                 \O_Y(\widetilde{D}_0' + 
              \sum_{i = 1}^d a_{i0} \widetilde{D}_i)
\]
and
\[
  \O_Y(k)\otimes 
  \O_Y(k_0      \widetilde{D}_0' + \ldots + k_s \widetilde{D}_s') = 
  \O_Y((k + k_0)\widetilde{D}_0' + \sum_{j = 1}^s k_j      
                \widetilde{D}_j' + \sum_{i = 1}^d k a_{i0} 
                \widetilde{D}_i).
\]
Each $\widetilde{m}\in\widetilde{M}_{\R}$ can be uniquely written as
$\widetilde{m} = m + \lambda_1\eta_1 + \ldots + \lambda_s\eta_s$, where
$m\in M_{\R}$, $\lambda_j\in\R$. The support polytope $\nabla$ is defined 
by the inequalities
\begin{eqnarray*}
  &&\<\widetilde{m},\widetilde{v}_i\> \ge -k a_{i0}, \quad i = 1,\ldots,d,\\
  &&\<\widetilde{m},\widetilde{e}_0\> \ge -k - k_0, \\
  &&\<\widetilde{m},\widetilde{e}_j\> \ge -k_j, \quad j = 1,\ldots,s
\end{eqnarray*}
in $\widetilde{M}_{\R}$. Let $\lambda_0 = k - \sum_{j = 1}^s \lambda_j$. 
Then the inequalities above are equivalent to 
\[
\{
  \<m,v_i\> \ge - \sum_{j = 0}^s \lambda_j a_{ij}, \quad i = 1,\ldots,d, \quad
  \lambda_j \ge -k_j, \quad j = 0,\ldots,s,\quad \sum_{j = 0}^s \lambda_j = k
\}.
\]
The same arguments as in \S 3 of \cite{Cattani-Cox-Dickenstein} 
show that $m = \sum_{j = 0}^s \lambda_j m_j$, where 
$m_j\in\Delta_j$, and $\widetilde{m}\in\nabla$ follows immediately.
\end{proof}
\begin{dfn} 
For any $J = \{j_1,\ldots,j_t\}\subset I := \{0,\ldots,s\}$
define the polytope
\[  
  \Delta_{j_1\ldots j_t} (k) := 
  \{
  \sum_{i = 1}^s \lambda_i \eta_i + 
  \sum_{t = 0}^s \lambda_t m_t:\, \lambda_j\ge 1,\, j\in J,\,
  \lambda_j\ge 0,\,j\in I\backslash J,\,\sum_{j = 0}^s \lambda_j = k,\,
  m_t\in\Delta_t 
  \}.
\] 
\end{dfn}
\begin{theo}For any sufficiently large integer $k > 0$
    \begin{eqnarray*}
     h^q(Y,\,\Omega_{Y/\P}^p(k)) 
      = \left\{
      \begin{array}{ll}
       \displaystyle
       \sum_{t = 0}^p (-1)^{p - t}
       \sum_{0 \le j_1 < \ldots < j_t \le s}
       \l (\Delta_{j_1 \ldots j_t}(k)), &q = 0,\\
       \\
       0,                               &q > 0.
      \end{array}
      \right.
    \end{eqnarray*}
\end{theo}
\begin{proof} From the $p$-th exterior power of the sequence (\ref{Omega_Y}),  
we have
\[
 0 \rightarrow \Omega^p_{Y/\P}(k) 
   \rightarrow \bigwedge\nolimits^p \pi^*(\E)\otimes \O_Y(k - p)
   \rightarrow \Omega^{p - 1}_{Y/\P}(k) 
   \rightarrow 0.
\]
Since $\O_Y(1)$ is ample, for all $k\gg 0$, we obtain the exact 
sequence of global sections
\[
 0 \rightarrow H^0(Y,\,\Omega^p_{Y/\P}(k)) 
   \rightarrow 
   H^0(Y,\,\bigwedge\nolimits^p \pi^*(\E)\otimes \O_Y(k - p))
   \rightarrow 
   H^0(Y,\,\Omega^{p - 1}_{Y/\P}(k)) 
   \rightarrow 0.
\]
Lemma 2.1 of \cite{Mavlyutov} implies the isomorphisms 
$\pi^*(\L_{j_m})\otimes\O_Y(-1)\simeq\O_Y(-\widetilde{D}_{j_m}')$
and consequently, the isomorphisms
\begin{eqnarray*}
  \bigwedge\nolimits^p \pi^*(\E)\otimes \O_Y(k - p) 
  &\simeq&
  \O_Y(k)\otimes\bigoplus_{0\le j_1 < \ldots < j_p \le s}
  \pi^*(\L_{j_1})\otimes\cdots\otimes
  \pi^*(\L_{j_p})\otimes\O_Y(-p) \\
  &\simeq& 
  \O_Y(k)\otimes\bigoplus_{0\le j_1 < \ldots < j_p \le s}
  \O_Y(-\widetilde{D}_{j_1}' - \ldots - \widetilde{D}_{j_p}').
\end{eqnarray*}
Hence
we obtain the relations
\[
  h^0(Y,\,\Omega^p_{Y/\P}(k)) = \sum_{0\le j_1 < \ldots < j_p \le s}
  \l(\Delta_{j_1\ldots j_p }(k)) - h^0(Y,\,\Omega^{p - 1}_{Y/\P}(k)).
\]
By induction we get the requested formula.
\end{proof}


\section{Appendix}

In the proof of Theorem \ref{comb_dual} we have used the following identity.

\begin{lem} For any face $F$ of a simple $n$-polytope $\Delta$
there are relations
\[ 
 \sum_{\Gamma \supseteq F} (-1)^{dim \Gamma} \binom{\dim\Gamma}{p}
 = 
 (-1)^{n} \binom{\dim F}{n - p}.
\] 
\end{lem}
\begin{proof} Since $\Delta$ is a simple, there are precisely 
\[
 \binom{n-s}{n-j}
\]
$j$-faces of $\Delta$ containing a given $s$-face $F$ of $\Delta$.
Hence, the requested formula is equivalent to the combinatorial identity
\[
 \sum_{j = 0}^n (-1)^j \binom{j}{p}\binom{n-s}{n-j} =
 (-1)^n \binom{s}{n-p},
\]
or
\[
 \sum_{j = 0}^{n-s} (-1)^j \binom{n-j}{p}\binom{n-s}{j} =
 \binom{s}{n-p}.
\]
We give a proof of this identity based on the method of integral 
representations. Rewrite the sum
\[
 S = \sum_{j = 0}^{n - s} (-1)^j \binom{n - j}{p}\binom{n - s}{j}
\]
as
\[
 S = \sum_{j = 0}^\infty (-1)^j \frac{1}{(2\pi i)^2}
 \iint_{\gamma} 
 \frac{(1 + z)^{n - j} (1 + w)^{n - s}}
 {z^{p + 1} w^{j + 1}}\, dz\, dw,
\]
where the cycle 
$\gamma$ is $\{|z| = \varepsilon\}\times \{|w| = \delta\}$. 
One can choose the numbers $\varepsilon > 0$ and $\delta > 0$ 
small enough, say $\varepsilon = 1/2$, $\delta = 3$, so
that the geometric series $\sum_{j = 0}^\infty (-(1+ z)^{-1}w^{-1})^{j}$
converges on the contour of integration $\gamma$, and it is possible 
to reverse the order of integration and summation. 
Therefore, summing up a geometric sequence,
\begin{eqnarray*} 
 S &=& \frac{1}{(2\pi i)^2} \iint_{\gamma} 
 \frac{(1 + z)^n (1 + w)^{n - s}}{z^{p + 1} w}
 \sum_{j = 0}^\infty 
 \left( \frac{-1}{(1 + z)w}\right)^j\, dz\, dw
 \\
 &=&
 \frac{1}{(2\pi i)^2}\iint_{\gamma}
 \frac{(1 + z)^n (1 + w)^{n - s}}{z^{p + 1} w 
 \left(1 +\displaystyle\frac{1}{(1 + z)w} \right)}\, dz\, dw
 \\
 &=&
 \frac{1}{2\pi i}\int_{|z| = 1/2} 
 \frac{(1 + z)^{n + 1}}{z^{p + 1}}\,dz \cdot
 \frac{1}{2\pi i}\int_{|w| = 3} 
 \frac{(1 + w)^{n - s}}{1 + zw + w}\,dw.
\end{eqnarray*}
By the residue theorem, the second integral is
\begin{eqnarray*}
 \frac{1}{2\pi i}\int_{|w| = 3}
 \frac{(1 + w)^{n - s}}{1 + zw + w}\,dw 
 &=&
 res_{w = w_0} \frac{(1 + w)^{n - s}}{(1 + z)(w + 1/(1 + z))}\\
 &=&
 \frac{\left(1 - \displaystyle\frac{1}{1 + z}\right)^{n - s}}{1 + z} =
 \frac{z^{n - s}}{(1 + z)^{n - s + 1}},
\end{eqnarray*}
where $w_0 = -1/(1 + z)\in\{|w|=3\}$, since 
$|w_0| \le \displaystyle\frac{1}{1 - 1/2} = 2 < 3$ for $|z| = 1/2$.
Finally,
\begin{eqnarray*}
 S &=& \frac{1}{2\pi i} \int_{|z| = 1/2}
 \frac{(1 + z)^{n + 1}}{z^{p + 1}} \cdot
 \frac{z^{n - s}}{(1 + z)^{n - s + 1}}\,dz \\
 &=&
 \frac{1}{2\pi i} \int_{|z|=1/2} 
 \frac{(1 + z)^s}{z^{p - n + s +1}} \,dz =
 \binom{s}{p -n + s} = \binom{s}{n - p},
\end{eqnarray*}
which concludes the proof.
\end{proof}


\begin{thebibliography}{99999}

\bibitem[Au]{Audin} M.~Audin, 
{\em The Topology of Torus Action on Symplectic
Manifolds}, Progress in Math. {\bf 93}, Birkh\"auser, Basel, Boston, Berlin,
1991.

\bibitem[BC]{Batyrev-Cox} V.V.~Batyrev, D.~Cox, {\em On the Hodge
structure of projective hypersurfaces in toric varieties}, Duke Math. J.
{\bf 75} (1995), 293--338.

\bibitem[Bo]{Bott} R.~Bott, {\em Homogeneous vector bundles},
Ann. of Math. {\bf 66} (1957), 203-248.

\bibitem[Br]{Broendsted} A.~Br\o ndsted, {\em An Introduction to Convex 
Polytopes}, Graduate Texts in Math., vol.~90, Springer, 1983.

\bibitem[CCD]{Cattani-Cox-Dickenstein} E.~Cattani, D.~Cox, A.~Dickenstein,
{\em Residues in toric varieties}, Composito Math. {\bf 108} (1997), 35-76.

\bibitem[C]{Cox} D.~Cox, {\em The homogeneous coordinate ring of a toric
variety}, J. Algebraic Geometry {\bf 4} (1995), 17-50.

\bibitem[D1]{Danilov} V.I.~Danilov, {\em The geometry of toric varieties},
Russian Math. Surveys {\bf 33} (1978), 97-154.

\bibitem[D2]{Danilov2} \bysame, {\em de~Rham complex on toroidal 
variety}, Algebraic geometry, (Chicago, IL, 1989), 26--38, 
Lecture Notes in Math., 1479, Springer, Berlin, 1991.

\bibitem[DK]{Danilov-Khovanskii} V.I.~Danilov, A.G.~Khovanski\^i,
{\em Newton polyhedra and an algorithm for computing
Hodge-Deligne numbers}, Math. USSR Izvestiya {\bf 29} (1987),
279-298.

\bibitem[De]{Dehn} M.~Dehn, {\em Die Eulersche Formel im Zusammenhang mit dem 
Inhalt in der nicht-Euklidiscen Geometrie}, 
Math. Ann. {\bf 61} (1905), 561-568.

\bibitem[Do]{Dolgachev} I.~Dolgachev, {\em Weighted projective varieties},
in {\it Group Action and Vector Fields}, ed. by J.B.~Carrell, Lecture Notes
in Math. {\bf 956}, Springer-Verlag, Berlin, 1982, 34-71.

\bibitem[EY]{Egorychev-Yuzhakov} G.P.~Egorychev, A.P.~Yuzhakov, {\em The
determination of generating functions and combinatorial sums by means of
multidimensional residues}, Siberian Math. J. {\bf 15} (1974).

\bibitem[E]{Ehrhart} E.~Ehrhart, {\em D\'emonstration de loi de 
r\'eciprocit\'e pour un poly\'edre entier}, C.R. Acad. Sci Paris
{\bf 265A} (1967), 5-7.

\bibitem[F1]{Fulton} W.~Fulton, {\em Intersection Theory}, Springer-Verlag,
1984.

\bibitem[F2]{Fulton.toric} \bysame, {\em Introduction to Toric Varieties},
Princeton Univ. Press, Princeton, NJ, 1993.

\bibitem[Ha]{Hartshorn} R.~Hartshorn, {\em Ample Subvarieties of Algebraic
Varieties}, Lecture Notes in Math. {\bf 156}, Springer-Verlag, Berlin,
Heidelberg, New York, 1970.

\bibitem[Hi]{Hirzebruch} F.~Hirzebruch, {\em Topological Methods in Algebraic
Geometry}, Springer-Verlag, Berlin, Heidelberg, New York, 1966.

\bibitem[Hu]{Huang} I-C.~Huang, {\em Cohomology of projective space seen by 
residual complex}, Trans. Amer. Math. Soc., to appear.

\bibitem[Is]{Ishida} M.-N.~Ishida, {\em Torus embeddings and dualizing
complex}, T\^ohoku Math. Jour. {\bf 32} (1980), 111-146.


\bibitem[K]{Kleiman} S.L.~Kleiman, {\em Toward a numerical theory of
ampleness}, Ann. of Math. {\bf 84} (1966), 293-344.

\bibitem[Mac]{Macdonald} I.G.~Macdonald, {\em Polynomials associated with
finite complexes}, J. London Math. Soc. (2) {\bf 4} (1971), 181-192.

\bibitem[Man]{Manin} Yu.I.~Manin, {\em Lectures on the $K$-functor in 
algebraic geometry}, Russ. Math. Surveys, {\bf 24(5)} (1969), 1-89.

\bibitem[MY]{Materov-Yuzhakov} E.N.~Materov, A.P.~Yuzhakov,
{\em The Bott formula for toric varieties and some combinatorial 
identies}, in {\em Complex analysis and differential operators 
(Russian)}, Krasnoyarsk, 2000, 85-92.

\bibitem[M1]{Mavlyutov} A.R.~Mavlyutov, {\em Cohomology of complete 
intersections in toric varieties}, Pacif. J. Math. {\bf 191} (1999), 133-144.

\bibitem[M2]{Mavlyutov2} \bysame, {\em On the chiral ring
of Calaby-Yau hypersurfaces in toric varieties}, preprint math.AG/0010318.

\bibitem[O1]{Oda.book} T.~Oda, {\em Convex Bodies and Algebraic Geometry --
an Introduction to the Theory of Toric Varieties}, Ergebnisse der Math. (3)
{\bf 15}, Springer-Verlag, Berlin, Heidelberg, New York, Paris, Tokyo, 1988.

\bibitem[O2]{Oda.de.Rham} \bysame, {\em The algebraic de Rham theorem for
toric varieties}, T\^ohoku Math. Jour. {\bf 45} (1993), 231-247.

\bibitem[OSS]{Okonek} C.~Okonek, M.~Schneider and H.~Spindler,
{\em Vector Bundles on Complex Projective Spaces}, Birkh\"auser, Boston,
Basel, Stuttgart, 1980.

\bibitem[Sn]{Snapper} E.~Snapper, {\em Multiples of divisors}, J. Math. and
Mech. {\bf 8} (1959), 967-922.

\bibitem[So]{Sommerville} D.M.Y.~Sommerville, {\em The relation connecting 
the angle-sums and volume of a polytope in space of $n$ dimensions},
Proc. Roy. Soc. London, Ser. A {\bf 115} (1927), 103-119.

\bibitem[St]{Stanley} R.~Stanley, {\em The number of the faces of a simplicial
convex polytope}, Advances in Math. {\bf 35} (1980), 236-238.

\bibitem[TE]{Torus.Embeddings} G.~Kempf, F.~Knudsen, D.~Mumford and
B.~Saint-Donat, {\em Toroidal Embeddings, I}, Lecture Notes in Math.
{\bf 339}, Springer-Verlag, Berlin, Heidelberg, New York, 1973.


\end{thebibliography}
\end{document}